\documentclass[12pt,reqno]{amsart}
\usepackage{amsmath, amsthm, amscd, amsfonts, amssymb, graphicx, color}
\usepackage[bookmarksnumbered, colorlinks, plainpages]{hyperref}

\theoremstyle{definition}

\theoremstyle{remark}

\numberwithin{equation}{section}
\begin{document}
\begin{center}
{\textbf{Yamabe Solitons on $(LCS)_{n}$-manifolds}}
\end{center}
\vskip 0.3cm
\begin{center}By\end{center}\vskip 0.3cm
\begin{center}
{Soumendu Roy \footnote{The first author is the corresponding author.},Santu Dey $^2$ and Arindam~~Bhattacharyya $^3$}
\end{center}
\vskip 0.3cm
\begin{center}
Department~of~Mathematics\\
Jadavpur~University,\\
Kolkata-700032,~India.\\
E-mail$^1$:soumendu1103mtma@gmail.com\\
E-mail$^2$: santu.mathju@gmail.com\\
E-mail$^3$: bhattachar1968@yahoo.co.in
\end{center}
\vskip 0.5cm
\begin{center}
\textbf{Abstract}\end{center}
The object of the present paper is to study some properties of $(LCS)_{n}$-manifolds whose metric is Yamabe soliton. We establish some characterization of $(LCS)_{n}$-manifolds when the soliton becomes steady. Next we have studied some certain curvature conditions of $(LCS)_{n}$-manifolds admitting Yamabe solitons. Lastly we construct a 3-dimensional $(LCS)_{n}$-manifold satisfying the results.
\vspace {0.3cm}
\section{\textbf{Introduction}}
In 2003, the notion of Lorentzian concircular structure manifolds (briefly,\\ $(LCS)_n$-manifolds) was introduced by Shaikh \cite{shaikh} with an example, that generalize the notion of LP-Sasakian manifolds which was introduced by Matsumoto \cite{matsumoto} and also by Mihai and Rosca \cite{mihai}. In 2005 and 2006, the application of $(LCS)_n$-manifolds to the general theory of relativity and cosmology was studied by Shaikh and Baishya \cite{baishya}, \cite{shba}.\\\\
Many other authors M. Ateceken\cite{atceken}, S. K. Hui\cite{hui}, D. Narain\cite{narain}, S. Yadav(\cite{yadav1}, \cite{yadav2}, \cite{yadav3}, \cite{yadav4}), A. Shaikh(\cite{basu}, \cite{basu1}) also studied the $(LCS)_n$-manifolds.\\\\
The concept of Yamabe flow was first introduced by Hamilton \cite{hamil} to construct Yamabe metrics on compact Riemannian manifolds. On a Riemannian or pseudo-Riemannian manifold $M$, a time-dependent metric $g(\cdot, t)$ is said to evolve by the Yamabe flow if the metric $g$ satisfies the given equation,
\begin{equation}\label{1.1}
  \frac{\partial }{\partial t}g(t)=-rg(t),\hspace{0.5cm} g(0)=g_{0},
\end{equation}
where $r$ is the scalar curvature of the manifold $M$.\\\\
In 2-dimension the Yamabe flow is equivalent to the Ricci flow (defined by $\frac{\partial }{\partial t}g(t) = -2S(g(t))$, where $S$ denotes the Ricci tensor). But in dimension $> 2$ the Yamabe and Ricci flows do not agree, since the Yamabe flow preserves the conformal class of the metric but the Ricci flow does not in general.\\\\
A Yamabe soliton correspond to self-similar solution of the Yamabe flow, is defined on a Riemannian or pseudo-Riemannian manifold $(M, g)$ by a vector field $\xi$ satisfying the equation \cite{barbosa},
\begin{equation}\label{1.2}
  \frac{1}{2}\pounds_\xi g = (r-\lambda)g,
\end{equation}
where $\pounds_\xi g$ denotes the Lie derivative of the metric $g$ along the vector field $\xi$, $r$ is the scalar curvature and $\lambda$ is a constant. Moreover a Yamabe soliton is said to be expanding if $\lambda < 0$, steady if $\lambda = 0$ and shrinking if $\lambda > 0$.\\
Yamabe solitons on a three-dimensional Sasakian manifold were studied by
R. Sharma \cite{sharma}.\\\\
Again,
\begin{equation}\label{1.3}
  R(X, Y )Z = \nabla_X\nabla_Y Z - \nabla_Y \nabla_X Z - \nabla_{[X,Y]} Z,
\end{equation}
\begin{eqnarray}\label{1.4}
 H(X, Y )Z &=& R(X, Y )Z -\frac{1}{(n-2)}[g(Y,Z)QX - g(X,Z)QY\nonumber \\
           &+& S(Y,Z)X - S(X,Z)Y],
\end{eqnarray}
 \begin{equation}\label{1.5}
       P(X, Y )Z = R(X, Y )Z -\frac{1}{(n-1)}[g(QY,Z)X - g(QX,Z)Y ],
     \end{equation}
 \begin{equation}\label{1.6}
   \tilde{C}(X,Y)Z=R(X,Y)Z-\frac{r}{n(n-1)}[g(Y,Z)X-g(X,Z)Y],
 \end{equation}
 \begin{equation}\label{1.7}
   W_2(X,Y)Z=R(X,Y)Z+\frac{1}{n-1}[g(X,Z)QY-g(Y,Z)QX],
 \end{equation}
 are the Riemannian-Christoffel curvature tensor $R$ \cite{ozgu}, the conharmonic
curvature tensor $H$ \cite{ishii}, the projective curvature tensor $P$ \cite{yano}, the concircular curvature tensor $\tilde{C}$ \cite{MISHRA}  and the $W_2$-curvature tensor \cite{MISHRA} respectively in a Riemannian manifold $(M^n, g)$, where $Q$ is the Ricci operator, defined by $S(X, Y) = g(QX, Y )$, $S$ is the Ricci tensor, $r = tr(S)$ is the scalar curvature, where $tr(S)$ is the trace of $S$ and $X, Y, Z \in \chi(M)$, $\chi(M)$ being the Lie algebra of vector fields of M.\\\\
In the present paper we study  Yamabe soliton on $(LCS)_{n}$-manifolds. The paper is organized as follows:\\
After introduction, section 2 is devoted for preliminaries on $(LCS)_{n}$-manifolds. In section 3, we have studied  Yamabe soliton on $(LCS)_{n}$-manifolds. Here we examined when $(LCS)_{n}$-manifold admits Yamabe soliton, then the manifold becomes $K-(LCS)_{n}$-manifold and Ricci symmetric. In this section we have also shown that, $(LCS)_{n}$-manifold admits Yamabe soliton is $\xi$-projectively flat, $\xi$-concircularly flat and $\xi$-conharmonically flat iff the soliton becomes steady. Section 4 deals with curvature properties on $(LCS)_{n}$-manifold admitting Yamabe soliton. Here we obtained some results on Yamabe soliton satisfying the conditions of the following type:\\
$S(\xi,X)\cdot R=0$, $S(\xi,X)\cdot W_2 = 0$, where $W_2$ is the $W_2$- curvature tensor and $S$ is the Ricci tensor. Also we have found that if the manifold admits Yamabe soliton then $R(\xi,X)\cdot S = 0$ and $W_2(\xi,X)\cdot S = 0$.\\
In the last section, we gave an example of 3-dimensional $(LCS)_{n}$-manifold on which we can easily verify our results.
\vspace {0.3cm}
\section{\textbf{Preliminaries}}
Let$(M,g)$ be an $n$-dimensional Lorentzian manifold admitting a unit timelike concircular vector field $\xi$. Then the vector field satisfying $g(\xi, \xi) = -1$. Since $\xi$ is a unit concircular vector field, it follows that there exists a nonzero 1-form $\eta$ such that, $g(X, \xi)= \eta(X)$. Also $\xi$ satisfies $\nabla \xi= \alpha(I + \eta\otimes\xi)$ with a nowhere zero smooth function $\alpha$ on $M$ verifying the equation $\nabla_{X} \alpha= (X\alpha)= d\alpha(X)=\rho\eta(X)$, for $\rho \in C^{\infty}(M)$, where $\nabla$ is the Levi-Civita connection of $g$ and $X$ is a vector field. Here also $\phi$ is the (1,1) tensor field, denoted by, $\phi:=\frac{1}{\alpha}\nabla\xi$.\\
The notion of Lorentzian para-Sasakian manifold was introduced by K. Matsumoto \cite{matsumoto}. More general the Lorentzian manifold M together with the unit timelike concircular vector field $\xi$, an  1-form $\eta$, and an (1,1) tensor field $\phi$ is said to be a Lorentzian concircular structure manifold $(M, g, \xi, \eta, \phi, \alpha)$ (briefly, $(LCS)_{n}$-manifold), which was introduced by A. A. Shaikh \cite{shaikh}.\\
In an $n$-dimensional $(LCS)_{n}$-manifold the following relations hold:\\
\begin{equation}\label{2.1}
  \phi^2=I+\eta\otimes\xi, \hspace{0.2cm} \eta(\xi)=-1, \hspace{0.2cm} \phi\xi=0, \hspace{0.2cm} \eta\circ\phi=0,
\end{equation}
\begin{equation}\label{2.2}
  g(\phi X,\phi Y)= g(X,Y)+\eta(X)\eta(Y) \hspace{0.2cm} and \hspace{0.2cm} g(\phi X, Y)=g(X,\phi Y),
\end{equation}
\begin{equation}\label{2.3}
  (\nabla_{X}\phi)Y=\alpha[g(X,Y)\xi+2\eta(X)\eta(Y)\xi+\eta(Y)X],
\end{equation}
for any $X,Y \in \chi(M)$.
\begin{equation}\label{new new}
  \phi X=\frac{1}{\alpha} \nabla_X \xi,
\end{equation}
\begin{equation}\label{2.4}
  \eta(\nabla_{X}\xi)=0, \hspace{0.2cm} \nabla_{\xi}\xi=0,
\end{equation}
\begin{equation}\label{2.5}
  R(X,Y)Z= (\alpha^2-\rho)[g(Y,Z)X-g(X,Z)Y],
\end{equation}
\begin{equation}\label{2.6}
  R(X,Y)\xi= (\alpha^2-\rho)[\eta(Y)X-\eta(X)Y],
\end{equation}
\begin{equation}\label{2.7}
  R(\xi,X)Y=(\alpha^2-\rho)[g(X,Y)\xi-\eta(Y)X],
\end{equation}
\begin{equation}\label{2.8}
  \eta(R(X,Y)Z)=(\alpha^2-\rho)[\eta(X)g(Y,Z)-\eta(Y)g(X,Z)],
\end{equation}
for any $X,Y,Z \in \chi(M)$.
\begin{equation}\label{2.9}
  \eta(R(X,Y)\xi)=0,
\end{equation}
\begin{equation}\label{new}
  S(X,Y)=(\alpha^2-\rho)(n-1)g(X,Y),
\end{equation}
\begin{equation}\label{new 1}
  r=n(n-1)(\alpha^2-\rho),
\end{equation}
\begin{equation}\label{2.10}
  \nabla\eta=\alpha(g+\eta\otimes\eta), \hspace{0.2cm} \nabla_{\xi}\eta=0,
\end{equation}
\begin{equation}\label{2.11}
  \pounds_\xi \phi=0, \hspace{0.2cm}\pounds_\xi \eta=0, \hspace{0.2cm}\pounds_\xi g=2\nabla\eta=2\alpha(g+\eta\otimes\eta),
\end{equation}
where $R$ is the Riemannian Curvature tensor, $S$ is the Ricci tensor, $r$ is the scalar curvature, $\nabla$ is the Levi-Civita  connection associated with $g$ and $\pounds_\xi$ denotes the Lie derivative along the vector field $\xi$.
\vspace {0.3cm}
\section{\textbf{Yamabe soliton on $(LCS)_n$ manifold}}
Let $(M, g, \xi, \eta, \phi, \alpha)$ be an $n$-dimensional $(LCS)_n$ manifold. Consider the Yamabe soliton on $M$ as:
\begin{equation}\label{3.1}
 \frac{1}{2}\pounds_\xi g = (r-\lambda)g.
\end{equation}
Then from \eqref{2.11}, we get,
\begin{equation}\label{3.2}
  \alpha[g(X,Y)+\eta(X)\eta(Y)]=(r-\lambda)g(X,Y),
\end{equation}
for all vector fields $X, Y$ on $M$.\\
which implies that,
\begin{equation}\label{3.3}
  (r-\lambda-\alpha)g(X,Y)-\alpha\eta(X)\eta(Y)=0.
\end{equation}
Taking $Y=\xi$ in the above equation and using \eqref{2.1}, we get,
\begin{equation}\label{3.4}
  (r-\lambda-\alpha)\eta(X)+\alpha\eta(X)=0.
\end{equation}
Then we have,
\begin{equation}\label{3.5}
  (r-\lambda)\eta(X)=0.
\end{equation}
Since $\eta(X) \neq 0$, so we get,
\begin{equation}\label{3.6}
  r=\lambda.
\end{equation}
Using the above equation in \eqref{3.1}, we have,
\begin{equation}\label{3.7}
  \pounds_\xi g=0.
\end{equation}
Thus $\xi$ is a Killing vector field and consequently $M$ is a $K-(LCS)_n$ manifold. Moreover, since $\lambda$ is constant, the scalar curvature $r$ is constant.\\\\
So we can state the following theorem:\\
\textbf{Theorem 3.1.} {\em If an $(LCS)_n$ manifold $(M, g, \xi, \eta, \phi, \alpha)$ admits a Yamabe soliton $(g, \xi)$, $\xi$ being the Reeb vector field of the Lorentzian concircular structure, then the scalar curvature is constant and the manifold is a $K-(LCS)_n$ manifold.}\\\\
Now from \eqref{new 1} and \eqref{3.6}, we get,
\begin{equation}\label{3.8}
  \lambda=n(n-1)(\alpha^2-\rho).
\end{equation}
Then using \eqref{new} and \eqref{3.8}, we obtain,
\begin{equation}\label{3.9}
  S(X,Y)=\frac{\lambda}{n}g(X,Y),
\end{equation}
for all vector fields $X, Y$ on $M$.\\\\
This leads to the following:\\
\textbf{Proposition 3.2.} {\em If an $(LCS)_n$ manifold $(M, g, \xi, \eta, \phi, \alpha)$ admits a Yamabe soliton $(g, \xi)$, $\xi$ being the Reeb vector field of the Lorentzian concircular structure, then the manifold becomes Einstein manifold.}\\\\
Now replacing the expression of $S$ from \eqref{3.9} in $$(\nabla_X S)(Y,Z) = X(S(Y,Z))-S(\nabla_X Y,Z)-S(Y,\nabla_X Z)$$ we get,
\begin{equation}\label{3.10}
  (\nabla_X S)(Y,Z) = \frac{\lambda}{n}(\nabla_X g)(Y,Z),
\end{equation}
which implies that,
\begin{equation}\label{3.11}
  \nabla S=0.
\end{equation}
This leads to the following:\\
\textbf{Proposition 3.3.} {\em If an $(LCS)_n$ manifold $(M, g, \xi, \eta, \phi, \alpha)$ admits a Yamabe soliton $(g, \xi)$, $\xi$ being the Reeb vector field of the Lorentzian concircular structure, then the manifold becomes Ricci symmetric.}\\\\
Again let the Ricci tensor $S$ of the $(LCS)_n$ manifold is $\eta$-recurrent i.e $$\nabla S = \eta \otimes S,$$ which implies that,
\begin{equation}\label{new 2}
  (\nabla_X S)(Y,Z)=\eta(X)S(Y,Z),
\end{equation}
for all vector fields $X, Y, Z$ on $M$.\\
Then using \eqref{3.11} and \eqref{3.9}, we get
\begin{equation}\label{new 3}
  \frac{\lambda}{n} \eta(x) g(Y,Z)=0.
\end{equation}
Taking $Y = \xi, Z = \xi$ in the above equation we obtain,
\begin{equation}\label{new 4}
  \lambda \eta(X)=0.
\end{equation}
As $\eta(X) \neq 0$, hence $\lambda = 0$.
Also from \eqref{3.6}, we get $r = 0$.\\\\
This leads to the following:\\
\textbf{Proposition 3.4.} {\em Let $(M, g, \xi, \eta, \phi, \alpha)$ be an $(LCS)_n$ manifold, admitting a Yamabe soliton $(g, \xi)$, $\xi$ being the Reeb vector field of the Lorentzian concircular structure. If the Ricci tensor $S$ of the manifold is $\eta$-recurrent (i.e $\nabla S = \eta \otimes S$), then the soliton is steady and the manifold becomes flat.}\\\\
Let us assume that a symmetric (0, 2) tensor field $h = \pounds_\xi g -2rg$ is parallel with respect to the Levi-Civita connection associated to $g$.\\
Then
\begin{equation}\label{3.12}
  h(\xi,\xi)=\pounds_\xi g(\xi,\xi) -2rg(\xi.\xi)=2\lambda,
\end{equation}
implies that,
\begin{equation}\label{3.13}
  \lambda=\frac{1}{2}h(\xi,\xi).
\end{equation}
Now as $h$ is parallel with respect to $g$, then from \cite{chandra}, we get,
\begin{equation}\label{3.14}
  h(X,Y)=-h(\xi,\xi)g(X,Y)
\end{equation}
for all vector fields $X, Y$ on $M$.\\
which leads to,
\begin{equation}\label{3.15}
  \pounds_\xi g(X,Y) = 2(r-\lambda)g(X,Y).
\end{equation}
So we can state the following theorem:\\
\textbf{Theorem 3.5.} {\em Let $(M, g, \xi, \eta, \phi, \alpha)$ be an $(LCS)_n$ manifold. Assume that a symmetric (0, 2) tensor field $h = \pounds_\xi g -2rg$ is parallel with respect to the Levi-Civita connection associated to $g$. Then $(g, \xi)$ yields an Yamabe soliton on $M$.}\\\\
We know,
 \begin{equation}\label{3.16}
    (\nabla_\xi Q)X = \nabla_\xi QX - Q(\nabla_\xi X),
 \end{equation}
 and
 \begin{equation}\label{3.17}
   (\nabla_\xi S)(X,Y) = \xi S(X,Y)-S(\nabla_\xi X,Y)-S(X,\nabla_\xi Y),
 \end{equation}
for any vector fields $X, Y$ on $M$.\\
Now using \eqref{3.9} we obtain,
\begin{equation}\label{3.18}
  QX=\frac{\lambda}{n}X,
\end{equation}
for any vector field $X$ on $M$.\\
Then in view of \eqref{3.9} and \eqref{3.18}, the equations \eqref{3.16} and \eqref{3.17} become
\begin{eqnarray}
  (\nabla_\xi Q)X &=& 0, \\
   (\nabla_\xi S)(X,Y)&=& 0,
\end{eqnarray}
respectively, for any vector fields $X, Y$ on $M$.\\\\
This leads to the following:\\
\textbf{Theorem 3.6.} {\em Let $(M, g, \xi, \eta, \phi, \alpha)$ be an $(LCS)_n$ manifold, admitting a Yamabe soliton $(g, \xi)$, $\xi$ being the Reeb vector field of the Lorentzian concircular structure. Then $Q$ and $S$ are parallel along $\xi$, where $Q$ is the Ricci operator, defined by $S(X, Y) = g(QX, Y )$ and $S$ is the Ricci tensor of $M$.}\\\\
Also in view of \eqref{3.18}, we obtain,
\begin{equation}\label{3.19}
  (\nabla_X Q)Y=\nabla_X QY - Q(\nabla_X Y)=0,
\end{equation}
for any vector fields $X, Y$ on $M$.\\\\
Then we have,\\
\textbf{Corollary 3.7.} {\em Let $(M, g, \xi, \eta, \phi, \alpha)$ be an $(LCS)_n$ manifold, admitting a Yamabe soliton $(g, \xi)$, $\xi$ being the Reeb vector field of the Lorentzian concircular structure. Then $Q$ is parallel to any arbitrary vector field on $M$.}\\\\
Let a Yamabe soliton is defined on an $n$-dimensional $(LCS)_n$ manifold $M$ as,
\begin{equation}\label{3.20}
  \frac{1}{2}\pounds_V g = (r-\lambda)g
\end{equation}
where $\pounds_V g$ denotes the Lie derivative of the metric $g$ along a vector field $V$ and $r, \lambda$ are defined as \eqref{1.2}.\\
Let $V$ be pointwise co-linear with $\xi$ i.e. $V = b\xi$ where $b$ is a function on $M$.\\ Then the equation \eqref{3.20} becomes,
\begin{equation}\label{3.21}
  \pounds_{b\xi} g(X,Y) = 2(r-\lambda)g(X,Y),
\end{equation}
for any vector fields $X, Y$ on $M$.\\
Applying the property of Lie derivative and Levi-Civita connection we have,
\begin{equation}\label{3.22}
  bg(\nabla_X \xi,Y)+(Xb)\eta(Y)+bg(\nabla_Y \xi,X)+(Yb)\eta(X)= 2(r-\lambda)g(X,Y).
\end{equation}
Using \eqref{new new}, the above equation reduces to,
\begin{equation}\label{3.23}
  b\alpha g(\phi X,Y)+(Xb)\eta(Y)+b\alpha g(\phi Y,X)+(Yb)\eta(X)= 2(r-\lambda)g(X,Y).
\end{equation}
Taking $Y = \xi$ in the above equation and using \eqref{2.1}, we get,
\begin{equation}\label{3.24}
  -Xb+(\xi b)\eta(X)=2(r-\lambda)\eta(X).
\end{equation}
Again putting $X = \xi$ in the above equation, we obtain,
\begin{equation}\label{3.25}
  \xi b=r-\lambda.
\end{equation}
Then using \eqref{3.25}, \eqref{3.24} becomes,
\begin{equation}\label{3.26}
  Xb=-(r-\lambda)\eta(X).
\end{equation}
Applying exterior differentiation in \eqref{3.26}, we have,
\begin{equation}\label{3.27}
  (r-\lambda)d\eta=0.
\end{equation}
Now in an $n$-dimensional $(LCS)_n$ manifold we have,\\ $$(d\eta)(X, Y) = X(\eta(Y)) - Y(\eta(X)) - \eta([X, Y]),$$ which implies
\begin{eqnarray}
  (d\eta)(X, Y) &=& g(Y,\nabla_X \xi) - g(X,\nabla_Y \xi) \nonumber \\
                &=& \alpha{g(Y, X) + \eta(X)\eta(Y)} - \alpha{g(X, Y) + \eta(X)\eta(Y)} \nonumber \\
                &=& 0.
\end{eqnarray}
Hence the 1-form $\eta$ is closed.\\
Then using the above equation, \eqref{3.27} implies that, either $r \neq \lambda$ or $r = \lambda.$\\
Now if $r \neq \lambda$ then from \eqref{3.20}, we have,
\begin{equation}\label{3.28}
  \pounds_V g = 2(r-\lambda)g,
\end{equation}
which implies $V$ is a conformal Killing vector field.\\
Again if $r = \lambda$ then from \eqref{3.26}, we get,
\begin{equation}\label{3.29}
  Xb=0,
\end{equation}
implies that $b$ is constant.\\\\
So we can state the following theorem:\\
\textbf{Theorem 3.8.} {\em Let $(M, g, \xi, \eta, \phi, \alpha)$ be an $(LCS)_n$ manifold, admitting a Yamabe soliton $(g, V)$, $V$ being a vector field on $M$. If $V$ is pointwise co-linear with $\xi$ then either $V$ is a conformal Killing vector field, provided $r \neq \lambda$, or $V$ is a constant multiple of $\xi$, where $\xi$ being the Reeb vector field of the Lorentzian concircular structure, $r$ is the scalar curvature and $\lambda$ is a constant. }\\\\
Also if $r = \lambda$ then from \eqref{3.20}, we obtain,
\begin{equation}\label{3.30}
  \pounds_V g=0,
\end{equation}
implies that $V$ is a Killing vector field.\\\\
Then we have,\\
\textbf{Corollary 3.9.} {\em Let $(M, g, \xi, \eta, \phi, \alpha)$ be an $(LCS)_n$ manifold, admitting a Yamabe soliton $(g, V)$, $V$ being a vector field on $M$. If $V$ is pointwise co-linear with $\xi$ and $r = \lambda$ then $V$ becomes a Killing vector field, where $\xi$ being the Reeb vector field of the Lorentzian concircular structure, $r$ is the scalar curvature and $\lambda$ is a constant.}\\\\
From the definition of projective curvature tensor \eqref{1.5}, defined on an \\$n$-dimensional $(LCS)_n$ manifold, we have,
\begin{equation}\label{3.31}
  P(X, Y )Z = R(X, Y )Z -\frac{1}{(n-1)}[S(Y,Z)X - S(X,Z)Y],
\end{equation}
for any vector fields $X, Y, Z$ on $M$.\\
Putting $Z = \xi$, we get,
\begin{equation}\label{3.32}
  P(X, Y )\xi = R(X, Y )\xi -\frac{1}{(n-1)}[S(Y,\xi)X - S(X,\xi)Y].
\end{equation}
Using \eqref{2.6} and \eqref{3.9}, we obtain,
\begin{equation}\label{3.33}
  P(X, Y )\xi= [(\alpha^2-\rho)-\frac{\lambda}{n(n-1)}][\eta(Y)X-\eta(X)Y].
\end{equation}
Again using \eqref{3.8}, we get,
\begin{equation}\label{3.34}
  P(X, Y )\xi=0.
\end{equation}
This leads to the following:\\
\textbf{Proposition 3.10.} {\em An $(LCS)_n$ manifold $(M, g, \xi, \eta, \phi, \alpha)$, admitting a Yamabe soliton $(g, \xi)$, $\xi$ being the Reeb vector field of the Lorentzian concircular structure, is $\xi$-projectively flat.}\\\\
From the definition of concircular curvature tensor \eqref{1.6}, defined on an \\$n$-dimensional $(LCS)_n$ manifold, we have,
\begin{equation}\label{3.35}
  \tilde{C}(X,Y)Z=R(X,Y)Z-\frac{r}{n(n-1)}[g(Y,Z)X-g(X,Z)Y],
\end{equation}
for any vector fields $X, Y, Z$ on $M$.\\
Putting $Z = \xi$, we get,
\begin{equation}\label{3.36}
  \tilde{C}(X,Y)\xi=R(X,Y)\xi-\frac{r}{n(n-1)}[g(Y,\xi)X-g(X,\xi)Y].
\end{equation}
Using \eqref{2.6} and \eqref{3.6}, we obtain,
\begin{equation}\label{3.37}
  \tilde{C}(X,Y)\xi=[(\alpha^2-\rho)-\frac{\lambda}{n(n-1)}][\eta(Y)X-\eta(X)Y].
\end{equation}
Again using \eqref{3.8}, we get,
\begin{equation}\label{3.38}
  \tilde{C}(X,Y)\xi=0.
\end{equation}
This leads to the following:\\
\textbf{Proposition 3.11.} {\em An $(LCS)_n$ manifold $(M, g, \xi, \eta, \phi, \alpha)$, admitting a Yamabe soliton $(g, \xi)$, $\xi$ being the Reeb vector field of the Lorentzian concircular structure, is $\xi$-concircularly flat.}\\\\
From the definition of conharmonic curvature tensor \eqref{1.4}, defined on an \\$n$-dimensional $(LCS)_n$ manifold, we have,
\begin{eqnarray}
  H(X, Y )Z &=& R(X, Y )Z -\frac{1}{(n-2)}[g(Y,Z)QX - g(X,Z)QY\nonumber \\
           &+& S(Y,Z)X - S(X,Z)Y],
\end{eqnarray}
for any vector fields $X, Y, Z$ on $M$.\\
Putting $Z = \xi$, we get,
\begin{eqnarray}
  H(X, Y )\xi &=& R(X, Y )\xi -\frac{1}{(n-2)}[g(Y,\xi)QX - g(X,\xi)QY\nonumber \\
           &+& S(Y,\xi)X - S(X,\xi)Y],
\end{eqnarray}
Using \eqref{2.6}, \eqref{3.9} and \eqref{3.18}, we obtain,
\begin{equation}\label{3.39}
  H(X, Y )\xi=[(\alpha^2-\rho)-\frac{2\lambda}{n(n-2)}][\eta(Y)X-\eta(X)Y].
\end{equation}
Again using \eqref{3.8}, we get,
\begin{equation}\label{3.40}
   H(X, Y )\xi=-\frac{\lambda}{(n-1)(n-2)}[\eta(Y)X-\eta(X)Y].
\end{equation}
This implies that $H(X, Y )\xi = 0$ iff $\lambda = 0$.\\\\
This leads to the following:\\
\textbf{Proposition 3.12.} {\em An $(LCS)_n$ manifold $(M, g, \xi, \eta, \phi, \alpha)$, admitting a Yamabe soliton $(g, \xi)$, $\xi$ being the Reeb vector field of the Lorentzian concircular structure, is $\xi$-conharmonically flat iff the soliton is steady.}\\
 \section{\textbf{Curvature properties on $(LCS)_n$ manifold admitting Yamabe soliton }}
We know,
\begin{equation}\label{4.1}
  R(\xi,X)\cdot S=S(R(\xi,X)Y,Z)+S(Y,R(\xi,X)Z),
\end{equation}
for any vector fields $X, Y, Z$ on $M$.\\
Using \eqref{2.7}, we obtain,
\begin{equation}\label{4.2}
   R(\xi,X)\cdot S=S((\alpha^2-\rho)(g(X,Y)\xi-\eta(Y)X),Z)+S(Y,(\alpha^2-\rho)(g(X,Z)\xi-\eta(Z)X)).
\end{equation}
Then using \eqref{3.9}, we get,
\begin{eqnarray}\label{4.3}
  R(\xi,X)\cdot S &=& \frac{\lambda}{n}(\alpha^2-\rho)[g(X,Y)\eta(Z)\nonumber\\
                  &-& g(X,Z)\eta(Y)+g(X.Z)\eta(Y)-g(X,Y)\eta(Z)],
\end{eqnarray}
which implies that, $$R(\xi,X)\cdot S = 0.$$
So we can state the following theorem:\\
\textbf{Theorem 4.1.} {\em If an $(LCS)_n$ manifold $(M, g, \xi, \eta, \phi, \alpha)$ admits a Yamabe soliton $(g, \xi)$, $\xi$ being the Reeb vector field of the Lorentzian concircular structure, then $R(\xi,X)\cdot S = 0$, i.e the manifold is $\xi$-semi symmetric, where $R$ is the Riemannian curvature tensor and $S$ is the Ricci tensor.}\\\\
Again the condition $S(\xi,X)\cdot R=0$ implies that,
\begin{multline}
  S(X,R(Y,Z)W)\xi-S(\xi,R(Y,Z)W)X+S(X,Y)R(\xi,Z)W-S(\xi,Y)R(X,Z)W \\
  +S(X,Z)R(Y,\xi)W-S(\xi,Z)R(Y,X)W+S(X,W)R(Y,Z)\xi-S(\xi,W)R(Y,Z)X\\
  =0,
\end{multline}
for any vector fields $X, Y, Z, W$ on $M$.\\
Taking the inner product with $\xi$, the above equation becomes,
\begin{multline}
  -S(X,R(Y,Z)W)-S(\xi,R(Y,Z)W)\eta(X)+S(X,Y)\eta(R(\xi,Z)W)\\
  -S(\xi,Y)\eta(R(X,Z)W)+S(X,Z)\eta(R(Y,\xi)W)-S(\xi,Z)\eta(R(Y,X)W)\\
  +S(X,W)\eta(R(Y,Z)\xi)-S(\xi,W)\eta(R(Y,Z)X)=0.
\end{multline}
Replacing the expression of $S$ from \eqref{3.9} and taking $Z=\xi$, $W=\xi$, we get,
\begin{multline}
\frac{\lambda}{n}[-g(X,R(Y,\xi)\xi)-\eta(R(Y,\xi)\xi)\eta(X)+g(X,Y)\eta(R(\xi,\xi)\xi)\\-\eta(Y)\eta(R(X,\xi)\xi)
+\eta(X)\eta(R(Y,\xi)\xi)-\eta(\xi)\eta(R(Y,X)\xi)+\eta(X)\eta(R(Y,\xi)\xi)\\
-\eta(\xi)\eta(R(Y,\xi)X)]=0.
\end{multline}
Now using \eqref{2.6}, \eqref{2.8}, \eqref{2.9}, we get on simplification,
\begin{equation}\label{4.4}
  \frac{\lambda}{n}(\alpha^2-\rho)[g(X,Y)+\eta(x)\eta(Y)]=0.
\end{equation}
Using \eqref{2.2}, the above equation becomes,
\begin{equation}\label{4.5}
  \frac{\lambda}{n}(\alpha^2-\rho)g(\phi X,\phi Y)=0,
\end{equation}
for any vector fields $X, Y$ on $M$.\\
This implies that,
\begin{equation}\label{4.6}
  \frac{\lambda}{n} (\alpha^2-\rho) = 0.
\end{equation}
Then using \eqref{3.8}, we get,
$$\frac{\lambda^2}{n^2(n-1)}=0,$$
implies that $\lambda = 0.$\\
Hence using \eqref{3.6}, $r = 0.$\\\\
So we can state the following theorem:\\
\textbf{Theorem 4.2.} {\em If an $(LCS)_n$ manifold $(M, g, \xi, \eta, \phi, \alpha)$ admits a Yamabe soliton $(g, \xi)$, $\xi$ being the Reeb vector field of the Lorentzian concircular structure and satisfies $S(\xi,X)\cdot R=0$ then the manifold becomes flat and the soliton is steady, where $R$ is the Riemannian curvature tensor and $S$ is the Ricci tensor.}\\\\
We know,
\begin{equation}\label{4.7}
  W_2(\xi,X)\cdot S=S(W_2(\xi,X)Y,Z)+S(Y,W_2(\xi,X)Z),
\end{equation}
for any vector fields $X, Y, Z$ on $M$.\\
Replacing the expression of S from \eqref{3.9} and using the definition of $W_2$-curvature tensor from \eqref{1.7}, we get,
\begin{eqnarray}\label{4.8}
   W_2(\xi,X)\cdot S &=& \frac{\lambda}{n}g(R(\xi,X)Y+\frac{1}{n-1}[g(\xi,Y)QX-g(X,Y)Q\xi],Z) \nonumber\\
                &+& \frac{\lambda}{n}g(Y,R(\xi,X)Z+\frac{1}{n-1}[g(\xi,Z)QX-g(X,Z)Q\xi]).\nonumber\\
\end{eqnarray}
Now using \eqref{2.7} and the property $g(QX, Y) = S(X ,Y)$, we obtain on simplification,
\begin{eqnarray}\label{4.9}
  W_2(\xi,X)\cdot S&=& \frac{\lambda}{n(n-1)}[\eta(Y)S(X,Z)-S(\xi,Z)g(X,Y)\nonumber\\
              &+& \eta(Z)S(X,Y)-S(\xi,Y)g(X,Z)].
\end{eqnarray}
Then using \eqref{3.9}, the above equation becomes,
\begin{equation}\label{4.10}
   W_2(\xi,X)\cdot S=0.
\end{equation}
So we can state the following theorem:\\
\textbf{Theorem 4.3.} {\em If an $(LCS)_n$ manifold $(M, g, \xi, \eta, \phi, \alpha)$ admits a Yamabe soliton $(g, \xi)$, $\xi$ being the Reeb vector field of the Lorentzian concircular structure, then $W_2(\xi,X)\cdot S = 0$, where $W_2$ is the $W_2$- curvature tensor and $S$ is the Ricci tensor.}\\\\
Again the condition $S(\xi,X)\cdot W_2 = 0$ implies that,
\begin{multline}
  S(X,W_2(Y,Z)V)\xi-S(\xi,W_2(Y,Z)V)X+S(X,Y)W_2(\xi,Z)V \\
  -S(\xi,Y)W_2(X,Z)V +S(X,Z)W_2(Y,\xi)V-S(\xi,Z)W_2(Y,X)V \\
  +S(X,V)W_2(Y,Z)\xi-S(\xi,V)W_2(Y,Z)X =0,
\end{multline}
for any vector fields $X, Y, Z, V$ on $M$.\\
Taking the inner product with $\xi$, the above equation becomes,
\begin{multline}
  -S(X,W_2(Y,Z)V)-S(\xi,W_2(Y,Z)V)\eta(X)+S(X,Y)\eta(W_2(\xi,Z)V) \\
  -S(\xi,Y)\eta(W_2(X,Z)V) +S(X,Z)\eta(W_2(Y,\xi)V)-S(\xi,Z)\eta(W_2(Y,X)V) \\
  +S(X,V)\eta(W_2(Y,Z)\xi)-S(\xi,V)\eta(W_2(Y,Z)X) =0.
\end{multline}
Replacing the expression of $S$ from \eqref{3.9} and taking $Z=\xi$, $V=\xi$, we get,
\begin{multline}
\frac{\lambda}{n}[-g(X,W_2(Y,\xi)\xi)-\eta(W_2(Y,\xi)\xi)\eta(X)+g(X,Y)\eta(W_2(\xi,\xi)\xi)\\
-\eta(Y)\eta(W_2(X,\xi)\xi)+\eta(X)\eta(W_2(Y,\xi)\xi)-\eta(\xi)\eta(W_2(Y,X)\xi)+\eta(X)\eta(W_2(Y,\xi)\xi)\\
-\eta(\xi)\eta(W_2(Y,\xi)X)]=0.
\end{multline}
Now using \eqref{1.7}, \eqref{2.6}, \eqref{2.8}, \eqref{new}, we obtain on simplification,
\begin{equation}\label{4.11}
  \frac{\lambda}{n}[g(X,Y)+\eta(X)\eta(Y)-(\alpha^2-\rho)g(X,Y)-(\alpha^2-\rho)\eta(X)\eta(Y)]=0,
\end{equation}
implies that,
\begin{equation}\label{4.12}
   \frac{\lambda}{n}(1-\alpha^2+\rho)[g(X,Y)+\eta(x)\eta(Y)]=0.
\end{equation}
Using \eqref{2.2}, the above equation becomes,
\begin{equation}\label{4.13}
  \frac{\lambda}{n}(1-\alpha^2+\rho)g(\phi X,\phi Y)=0.
\end{equation}
for any vector fields $X, Y$ on $M$.\\
This implies that, $$\lambda (1-\alpha^2+\rho)=0.$$
Then either $\lambda = 0,$ or $\alpha^2 - \rho = 1.$\\
Now if $\alpha^2 - \rho = 1$, then from \eqref{new 1}, we have,$$r=n(n-1).$$
So we can state the following theorem:\\
\textbf{Theorem 4.4.} {\em If an $(LCS)_n$ manifold $(M, g, \xi, \eta, \phi, \alpha)$ admits a Yamabe soliton $(g, \xi)$, $\xi$ being the Reeb vector field of the Lorentzian concircular structure and satisfies $S(\xi,X)\cdot W_2 = 0$ then either the soliton is steady, or $r = n(n-1)$, where $W_2$ is the $W_2$- curvature tensor and $r$ is the scalar curvature.}\\\\
\section{\textbf{Example of an $(LCS)_3$  manifold:}}
We consider the 3-dimensional manifold $M = \{(x, y, z) \in \mathbb{R}^3, z \neq 0 \}$,
where $(x, y, z)$ are standard coordinates in $\mathbb{R}^3$. Let ${e_1, e_2, e_3}$ be a linearly independent system of vector fields on $M$ given by,
\begin{equation}
   e_1=z\frac{\partial}{\partial x},\quad e_2 =z \frac{\partial}{\partial y}, \quad e_3 =z\frac{\partial}{\partial z}. \nonumber
\end{equation}
Let $g$ be the Riemannian metric defined by,
\begin{equation}\label{5.1}
  g(e_1,e_1)= g(e_2,e_2) = 1, \quad g(e_3,e_3)=-1,\nonumber
\end{equation}
\begin{equation}\label{5.2}
   g(e_1,e_2) = g(e_2,e_3)= g(e_3,e_1) =0.\nonumber
\end{equation}
Let $\eta$ be the 1-form defined by $\eta(Z) = g(Z,e_3)$, for any $Z \in \chi(M)$,where $\chi(M)$ is the set of all differentiable vector fields on $M$ and $\phi$ be the (1, 1)-tensor field defined by,
\begin{equation}
  \phi e_1=e_2, \quad \phi e_2=e_1,\quad  \phi e_3=0.\nonumber
\end{equation}
Then, using the linearity of $\phi$ and $g$, we have $$\eta(e_3) = -1, \phi ^2 (Z) = Z + \eta(Z)e_3$$ and $$g(\phi Z,\phi W) = g(Z,W) + \eta(Z)\eta(W),$$ for any $Z,W \in \chi(M)$.\\
Let $\nabla$ be the Levi-Civita connection with respect to the Riemannian metric $g$. Then we have,
  $$ [e_1,e_2] = 0, [e_2,e_3] = -e_2, [e_1,e_3] = -e_1.$$
The connection $\nabla$ of the metric $g$ is given by,
\begin{eqnarray}
  2g(\nabla_X Y,Z) &=& Xg(Y,Z)+Yg(Z,X)-Zg(X,Y)\nonumber \\
                   &-& g(X, [Y,Z])-g(Y, [X, Z]) + g(Z, [X, Y]),\nonumber
\end{eqnarray}
which is known as Koszul’s formula.\\
Using Koszul’s formula, we can easily calculate,
$$\nabla_{e_1} e_3 = -e_1,  \nabla_{e_2} e_3 = -e_2 ,  \nabla_{e_3} e_3 = 0 ,$$
$$\nabla_{e_1} e_1 = -e_3,  \nabla_{e_2} e_1 = 0,  \nabla_{e_3} e_1 = 0,$$
$$\nabla_{e_1} e_2 = 0,  \nabla_{e_2} e_2 = -e_3,  \nabla_{e_3} e_2 = 0.$$
Hence in this case $(g, \xi, \eta, \phi, \alpha)$ is an $(LCS)_3$-structure on $M$, where $\alpha = -1$.\\
Also as $\alpha = -1$ then $\rho = 0$ and consequently, $(M, g, \xi, \eta, \phi, \alpha)$ is an $(LCS)_n$ manifold of dimension $3$.\\
On this manifold $(M, g, \xi, \eta, \phi, \alpha)$, we can easily verify our results.\\\\

\end{document}